\documentstyle[amscd,amssymb,verbatim,11pt]{amsart}

\theoremstyle{plain}
\newtheorem{Prop}{Proposition}[section]
\newtheorem{Thm}[Prop]{Theorem}
\newtheorem{Cor}[Prop]{Corollary}

\newtheorem{Lem}[Prop]{Lemma}

\theoremstyle{definition}
\newtheorem{Def}[Prop]{Definition}

\theoremstyle{remark}

\def\asdim{\text{asdim}}
\def\diam{\text{diam}}
\def\ccc{{\mathfrak{c}}}
\def\CC{{\mathfrak{C}}}
\def\CCoP{{{\mathfrak{C}}^\circ_P}}
\def\CCoQ{{{\mathfrak{C}}^\circ_Q}}

\def\RR{{\mathbb R}}

\def\PP{{\mathcal P}}
\def\UU{{\mathcal U}}

\def\ZZ{{\mathbb Z}}

\begin{document}
\title[Assouad-Nagata dimension of tree-graded spaces]{Assouad-Nagata dimension of tree-graded spaces}

\author{N.~Brodskiy}
\address{University of Tennessee, Knoxville, TN 37996, USA}
\email{brodskiy@@math.utk.edu}

\author{J.~Higes}
\address{Departamento de Geometr\'{\i}a y Topolog\'{\i}a,
Facultad de CC.Matem\'aticas. Universidad Complutense de Madrid. Madrid, 28040 Spain}
\email{josemhiges@@yahoo.es}



\thanks{ The second named author is supported by project MEC, MTM2006-0825.}

\begin{abstract}
Given a metric space $X$ of finite asymptotic dimension $\asdim X\le n$,
we consider a quasi-isometric invariant of the space called dimension function.
The space is said to have asymptotic Assouad-Nagata dimension $\asdim_{AN} X \le n$
if there is a linear dimension function in this dimension.
We prove that if $X$ is a tree-graded space (as introduced by C.~Drutu and M.~Sapir)
and for some positive integer $n$ a function $f$ serves as
an $n$-dimensional dimension function for all pieces of $X$, then the function
$300\cdot f$ serves as an $n$-dimensional dimension function for $X$.
As a corollary we find a formula for the asymptotic Assouad-Nagata dimension of the free product
of finitely generated infinite groups: $$\asdim_{AN} (G*H)=\max\{\asdim_{AN} (G), \asdim_{AN} (H)\}.$$
\end{abstract}

\maketitle


\section{Introduction and Preliminaries}

Asymptotic dimension was introduced by Gromov in~\cite{Gro asym invar} as a large scale invariant of a metric space. Any finitely generated group can be equipped with a word metric. The idea of Gromov was that asymptotic dimension is an invariant of the finitely generated group, i.e. does not depend on the word metric. The linear version of asymptotic dimension was called asymptotic dimension of linear type~\cite{Roe lectures}, asymptotic dimension with Higson property~\cite{Dran-Zar} and more recently asymptotic Assouad-Nagata dimension.

Spaces of finite asymptotic Assouad-Nagata dimension have some extra properties that spaces of finite asymptotic dimension do not necessarily have. For example if a metric space is of finite asymptotic Assouad-Nagata dimension then it satisfies nice Lipschitz extension properties (see~\cite{Lang-Sch Nagata dim} and~\cite{Brod-Dydak-Higes-Mitra}). It was proved in~\cite{Dydak-Higes} that the asymptotic Assouad-Nagata dimension bounds the topological dimension of every asymptotic cone of a metric space.  As a consequence of the results of Gal in~\cite{Gal} every metric space of finite asymptotic Assouad-Nagata dimension has Hilbert space compression one.

The list of finitely generated groups of finite asymptotic Assouad-Nagata dimension is somewhat limited. The classes of groups known to be of finite asymptotic Assouad-Nagata dimension are Coxeter groups, abelian groups, hyperbolic groups, free groups, and some types of Baumslag-Solitar groups. It is unknown if nilpotent groups are of finite asymptotic Assouad-Nagata dimension. 

A natural problem is to study the behavior of asymptotic Assouad-Nagata dimension under such operations on groups as free product, HNN-extension, and free product with amalgamation. This problem was solved in~\cite{BD} in case of asymptotic dimension. A formula for the asymptotic dimension of the free product of two groups was introduced in~\cite{DranFree} by Bell, Dranishnikov and Keesling. Unfortunately the techniques used in that paper cannot be extended directly to asymptotic Assouad-Nagata dimension.

The goal of the present paper is to obtain a formula for the asymptotic Assouad-Nagata dimension of the free product of two groups. This result will answer~\cite[Question 8.8]{Brod-Dydak-Levin-Mitra}. The techniques used in our proof allow us to extend the results to a larger class of spaces, the so called Tree-graded spaces.

\subsection{Dimension}

Let $(X,d)$ be a metric space and let $r$ be a positive real number. An {\it $r$-scale chain} (or $r$-chain) between two points $x$ and $y$ is defined as a finite sequence of points $\{x= x_0, x_1, ..., x_m = y\}$ such that $d(x_i, x_{i+1}) < r$ for every $i = 0, ..., m-1$.
A subset $S$ of a metric space $(X,d)$ is said to be {\it $r$-scale connected} if every two elements of $S$ can be connected by an $r$-scale chain
contained in $S$. 

Let $\UU$ be a subset of a metric space $(X,d)$. 
A maximal $r$-scale connected subset of $\UU$ is called an {\it $r$-scale connected component} of $\UU$.
Notice that any set $\UU$ is a union of its $r$-scale connected components and the distance between any two of the components is at least $r$.

\begin{Def}
Let $s$ and $M$ be two positive numbers. A metric space $(X,d)$ is said to have {\it dimension $\le n$ at the scale $s$ of magnitude $M$} if there is a cover $\UU =\{\UU_0, ...,\UU_n\}$ of $X$ such that every $s$-scale connected component of each $\UU_i$ is $M$-bounded (i.e. the diameter of every component is bounded by $M$).
\end{Def}

An increasing function $f_X: \RR_+ \to \RR_+$ is said to be an {\it $n$-dimensional control function} of $(X, d)$ if for any $s>0$ the space $(X,d)$ has dimension $\le n$ at the scale $s$ of magnitude $f_X(s)$.

A metric space $(X,d)$ is said to be of {\it asymptotic dimension} at most $n$ (notation: $\asdim (X, d) \le n$) if there exists an  $n$-dimensional control function for $(X,d)$. 

\begin{Def}
A metric space $(X,d)$ is said to be of {\it asymptotic Assouad-Nagata dimension} at most $n$ (notation: $\asdim_{AN}(X, d) \le n$) if there exists a linear $n$-dimensional control function i.e. there exists two positive constants $C$ and $b$ such that the function  $f(x) = C \cdot x+b$ is an $n$-dimensional control function for $(X,d)$.
\end{Def}

\subsection{Tree-graded spaces}

Tree-graded spaces were introduced by C.~Drutu and M.~Sapir~\cite{Drutu-Sapir_Tree-graded spaces} in order to study relatively hyperbolic groups.
The main interest in the notion of tree-graded space resides in the following characterization of relatively
hyperbolic groups proved by C.~Drutu, D.~Osin, and M.~Sapir: A finitely generated group $G$ is relatively hyperbolic with respect to finitely
generated subgroups $H_1, \dots , H_n$ if and only if every asymptotic cone $Con^\omega(G; e, d)$ is tree-graded with
respect to $\omega$-limits of sequences of cosets of the subgroups $H_i$~\cite{Drutu-Sapir_Tree-graded spaces}.

Tree-graded spaces are generalizations of $\RR$-–trees and have many properties similar to the properties of $\RR$–-trees.
They appear naturally as asymptotic cones of groups.
A theory of actions of groups on tree-graded spaces is developed
recently by C.~Drutu and M.~Sapir~\cite{Drutu-Sapir_Groups acting on tree-graded spaces}.

\begin{Def}\label{def tree-graded}
Let $X$ be a complete geodesic metric space and let $\PP$ be a collection
of closed geodesic subsets (called {\it pieces}) covering $X$.
The space $X$ is called {\it tree-graded with respect to $\PP$} if the following two properties are satisfied:
\begin{itemize}
\item[$(T_1)$] Every two different pieces have at most one common point.
\item[$(T''_2)$] Every simple loop in $X$ is contained in one piece.
\end{itemize}
\end{Def}

The main result of the paper is the following theorem.

\begin{Thm}
Let $X$ be a geodesic tree-graded space.
Let $n$ be a positive integer and $r$ and $f(r)\ge r$ be positive real numbers.
If each piece of $X$ has dimension $\le n$ at the scale $r$ of magnitude $f(r)$,
then the space $X$ has dimension $\le n$ at the scale $r$ of magnitude $300\cdot f(r)$.
\end{Thm}

A simple example of a tree-graded space is a Cayley graph of a free product $G*H$ of finitely generated groups when
the generating set for $G*H$ is taken as the union of generating sets for $G$ and $H$.

\begin{Cor}
Let $G$ and $H$ be finitely generated infinite groups. Then 
$$\asdim_{AN} (G*H)=\max\{\asdim_{AN} (G), \asdim_{AN} (H)\}.$$
\end{Cor}

Notice that if two groups $G$ and $H$ are finite, the free product $G*H$ is quasi-isometric to a tree, thus has asymptotic Assouad-Nagata dimension 1.
If only one of the groups is finite, say $G$, then we can reduce the computation of the asymptotic Assouad-Nagata dimension to the infinite case as follows: 
$asdim_{AN}(H)\le asdim_{AN}( G*H) \le asdim_{AN}((G\times \ZZ)*H)\le asdim_{AN}(H)$. So, the only interesting case is when both groups are infinite. 

\subsection{Idea of the proof of the main theorem}

Let $X$ be a geodesic metric space which is tree-graded with respect to a collection $\PP=\{P_\lambda\}_{\lambda\in\Lambda}$ of pieces.
Fix a base point $\widehat{x}\in X$. For every piece $P$ we let the projection of $\widehat{x}$ to $P$ be the base point of $P$ (call it $x_P$).
Since each piece $P$ has dimension $\le n$ at the scale $r$ of magnitude $f(r)$, we fix a coloring $c_P\colon P\to \{0,1,2,\dots,n\}$ of $P$ using $n+1$ colors such that the base point $x_P$ has color~$0$.

Let us try now to produce a coloring $c_X\colon X\to \{0,1,2,\dots,n\}$ of $X$.
Put $c_X(\widehat{x})=0$. Given any point $x\in X$ consider a geodesic $\gamma$ from $\widehat{x}$ to $x$.
As this geodesic travels from $\widehat{x}$ to $x$ it passes through different pieces of $X$ by entering a piece $P$ at the base point $x_P$ and exiting $P$ at (possibly) different point $\pi_P(x)$ which is the projection of the point $x$ to the piece $P$.
The first idea is to add up all the changes of colors along the way $\gamma$ (modulo $n+1$):
$$ c_X(x)=\sum_P c_P(\pi_P(x)) \quad (\text{mod } n+1)$$
Notice that this sum becomes finite if we require each coloring $c_P$ to be zero on the ball $B(x_P,2r)$.
The problem with this coloring $c_X$ is that the path $\gamma$ may never change color
if it changes pieces very often (say stays in one piece for no longer than $r$).
Now we improve the formula by introducing additional changes of color: if the path $\gamma$ does not change color for longer than $99\cdot f(r)$, we make it change the color. Namely, let $\beta_1,\dots, \beta_k$ be (maximal) subintervals of the path $\gamma$ such that each $\beta_i$ has its interior $c_X$-colored in one color and the length $|\beta_i|$ of $\beta_i$ is at least $99\cdot f(r)$.
$$ c'_X(x)=\sum_P c_P(\pi_P(x))+\sum_{i=1}^k \left\lfloor \frac{|\beta_i|}{99\cdot f(r)}\right\rfloor \quad (\text{mod } n+1)$$
The reason this coloring may not work is that its restriction to a piece $P$ of the space $X$ may differ from the coloring $c_P$.
It may happen because the second sum in the formula for $c'_X$ has nothing to do with the tree-graded structure of $X$.
Our last modification of the formula for the coloring uses the following definition.

\begin{Def}\label{notation floor}
Let $\beta$ be any directed geodesic in a tree-graded space $X$ from a point $\beta_-$ to a point $\beta_+$.
Consider the induced tree-graded structure on $\beta$. Let $[\beta',\beta_+]$ be the piece of this structure containing $\beta_+$.
We denote by $\lfloor\beta\rfloor$ the subpath of $\beta$ from $\beta_-$ to $\beta'$.
\end{Def}

Here is the coloring proving that the space $X$ has dimension $\le n$ at the scale $r$ of magnitude $100\cdot f(r)$.

$$ c''_X(x)=\sum_P c_P(\pi_P(x))+\sum_{i=1}^k \left\lfloor \frac{|\lfloor\beta_i\rfloor|}{99\cdot f(r)}\right\rfloor \quad (\text{mod } n+1)$$

We will give a slightly different and direct definition of a coloring of $X$ in Section~\ref{section Proof of the main theorem} and will prove the theorem there.

\section{Tree-graded spaces}\label{section Tree-graded spaces}

We collect in this section all the properties of tree-graded spaces that we use in the proof of our main theorem.

\begin{Lem}[{\cite[Proposition 2.17]{Drutu-Sapir_Tree-graded spaces}}]\label{lemma geodesic not missing a cut point}
Condition $(T''_2)$ in the definition of tree-graded spaces can be replaced by the following condition:
\begin{itemize}
\item[$(T'_2)$] For every topological arc $C \colon [0, d] \to X$ and $t\in [0, d]$, let $C[t-a, t + b]$ be a maximal sub-arc
of $C$ containing $C(t)$ and contained in one piece. Then every other topological arc with the same
endpoints as $C$ must contain the points $C(t - a)$ and $C(t + b)$.
\end{itemize}
\end{Lem}

If $X$ is tree-graded with respect to $\PP$ then we can always
add some or all one-point subsets of $X$ to $\PP$, and $X$ will be tree-graded with respect to a bigger set of pieces.
To avoid using extra pieces, we shall assume that pieces cannot contain other pieces.

Given a geodesic $\gamma$ in $X$ and a piece $P$ of $X$, we denote by $\gamma_P$ the intersection of $\gamma$ with $P$.
By~\cite[Corollary 2.10]{Drutu-Sapir_Tree-graded spaces}, $\gamma_P$ is either empty or a point or a closed subpath of $\gamma$.
Notice that $\gamma_P$ is a piece of the tree-graded structure induced on $\gamma$ by the structure on $X$.

\begin{Lem}[{\cite[Lemma 2.6]{Drutu-Sapir_Tree-graded spaces}}]\label{lemma d(x,P)}
Let $X$ be a tree-graded space and $P$ be a piece of $X$.
For every point $x\in X$ there exists a unique point $y\in P$ such that $d(x,P) = d(x,y)$.
Moreover, every geodesic joining $x$ with a point of $P$ contains $y$.
\end{Lem}

\begin{Def}
The point y in~\ref{lemma d(x,P)} is called {\it the projection of x onto the piece P} and is denoted $\pi_P(x)$.
\end{Def}

We will use the following lemma several times in the following context:
if $d(x,y)\le r$ and $\pi_P(x)\ne \pi_P(y)$, then $d(x,P)\le r$.

\begin{Lem}[{\cite[Lemma 2.8]{Drutu-Sapir_Tree-graded spaces}}]\label{lemma two points project to one}
Let $X$ be a tree-graded space and $P$ be a piece of $X$.
If $x$ and $y$ are two points of $X$ with $d(x,y)\le d(x,P)$, then $\pi_P(x)=\pi_P(y)$.
\end{Lem}

\begin{Lem}[{\cite[Lemma 2.20]{Drutu-Sapir_Tree-graded spaces}}]\label{lemma retraction to piece}
Let $X$ be a tree-graded space and $P$ be a piece of $X$.
The projection $\pi_P\colon X\to P$ is a 1-Lipschitz retraction.
\end{Lem}



\begin{Lem}\label{lemma chain not missing cut point}
Let $X$ be a tree-graded space and $\{x_i\}_{i=0}^m$ be an $r$-chain of points in $X$.
Let $\omega$ be a geodesic from $x_0$ to $x_m$.
If $P$ is a piece of $X$ such that $\omega_P\ne\emptyset$,
then $\omega_P$ is a path from $\omega_P^-=\pi_P(x_0)$ to $\omega_P^+=\pi_P(x_m)$ and
there is a point $x_k$ of the chain such that
$\pi_P(x_k)=\omega_P^-$ and $d(x_k,\omega_P^-)\le r$.
Moreover, any geodesic from $x_k$ to $x_m$ passes through $\omega_P^-$.
\end{Lem}

\begin{pf}
By~\ref{lemma geodesic not missing a cut point}, any geodesic from $x_0$ to $x_m$ passes through $\omega_P^-$.
We consider the maximal index $k$ such that any geodesic from $x_k$ to $x_m$ passes through $\omega_P^-$ and $\pi_P(x_k)=\omega_P^-$.
If $\pi_P(x_{k+1})\ne \omega_P^-$, then $d(x_k,P)\le r$ by~\ref{lemma two points project to one};
thus $d(x_k,\omega_P^-)=d(x_k,\pi_P(x_k))=d(x_k,P)\le r$.
If there is a geodesic from $x_{k+1}$ to $x_m$ not passing through $\omega_P^-$,
then by~\ref{lemma geodesic not missing a cut point}, any geodesic from $x_k$ to $x_{k+1}$ passes through $\omega_P^-$.
Since $d(x_k,x_{k+1})\le r$, we have $d(x_k,\omega_P^-)\le r$.
\end{pf}

\section{Proof of the main theorem}\label{section Proof of the main theorem}

Let us introduce some notations before we define a coloring $c^*_X\colon X\to \{0,1,2,\dots,n\}$ of $X$.
Fix a base point $\widehat{x}\in X$. For every piece $P$ we let the projection of $\widehat{x}$ to $P$ be the base point of $P$ (call it $x_P$).

Since each piece $P$ has dimension $\le n$ at the scale $r$ of magnitude $f(r)$,
we fix a coloring $c'_P\colon P\to \{0,1,2,\dots,n\}$ of $P$ using $n+1$ colors
such that $r$-components of each color are $f(r)$-bounded.
Denote by $c_P$ the coloring of $P$ obtained from $c'_P$ by changing the color of
the closed $2r$-neighborhood $B(x_P,2r)$ of the base point $x_P$ to the color 0.
Notice that $r$-components of any $c_P$-color other than 0 are $f(r)$-bounded and
$r$-components of $c_P$-color 0 are $8\cdot f(r)$-bounded.
Since the component of the point $x_P$ of $c_P$-color 0 will play an important role
in our definition of the coloring $c^*_X$, we give it a special notation $\CC_P$.
Denote by $\CCoP$ the subset $\CC_P\setminus\{x_P\}$.

Given any point $x\in X$ consider a geodesic $\gamma$ from $\widehat{x}$ to $x$.
The intersection $\gamma_P$ of $\gamma$ with a piece $P$ is either empty or a closed subpath of $\gamma$ from $x_P$ to $\pi_P(x)$.

\begin{Def}
We call $\gamma_P$ a {\it short} piece of $\gamma$ if its endpoint $\pi_P(x)$ belongs to the set $\CC_P$.
Otherwise we call this piece of $\gamma$ {\it long}.
\end{Def}

For every long piece $\gamma_P$ we call the {\it reduced long piece} to the piece  $\gamma_P' = \gamma_P \setminus B(x_P, \frac{r}{2})$.
Notice that any short piece $\gamma_P$ has length at most $8\cdot f(r)$, any long piece has length at least $2r$ and any reduced long piece has length at least $\frac{3\cdot r}{2}$.

Consider the complement in $\gamma$ of the interiors of all reduced long pieces of $\gamma$.
Since there are only finitely many long pieces of $\gamma$, this complement is a union of closed subpaths $\{\beta_i\}_{i=1}^k$ of $\gamma$.
Using Definition~\ref{notation floor}, we define a coloring of $X$ as follows:
$$ c^*_X(x)=\sum_P c_P(\pi_P(x))+\sum_{i=1}^k \left\lfloor \frac{|\lfloor\beta_i\rfloor|}{99\cdot f(r)}\right\rfloor \quad (\text{mod } n+1)$$

It is easy to check using~\ref{lemma geodesic not missing a cut point} that the color $c^*_X(x)$
does not depend on the choice of the geodesic $\gamma$ from $\widehat{x}$ to $x$.
The following Lemma claims that (except possibly at the base point $x_P$)
the coloring $c^*_X$ differs from the coloring $c_P$ by a constant modulo $n+1$.

\begin{Lem}\label{lemma c*X almost equals cP}
For any piece $P$ of the tree-graded space $X$ we have either
$$ c^*_X|_{P\setminus x_P} = c_P + c^*_X(x_P) \quad (\text{mod } n+1)$$
or
$$ c^*_X|_{P\setminus x_P} = c_P + c^*_X(x_P) + 1 \quad (\text{mod } n+1)$$
In particular, in any piece $P$ all $r$-components of each $c^*_X$-color are $8\cdot f(r)$-bounded.
\end{Lem}

\begin{pf}
Let us fix a piece $Q$ of $X$, a point $x\in Q\setminus\{x_Q\}$, and a geodesic $\gamma$ from $\widehat{x}$ to $x$.
Clearly $\sum_P c_P(\pi_P(x))=c_Q(x)+\sum_{P\ne Q} c_P(\pi_P(x))$ and
the first sum in the definition of the coloring $c^*_X$ changes by $c_Q(x)$ as we go from $x_Q$ to $x$. Let $Q-1$ denote the piece that contains $x_Q$ and a small portion of $\gamma$ just before $X_Q$.
If the piece $\gamma_{Q-1}$ is long, then $ c^*_X(x) = c_Q(x) + c^*_X(x_Q) \quad (\text{mod } n+1)$ because the second sum in the definition of the coloring $c^*_X$ has the same value for the points $x$ and $x_Q$.

If the piece $\gamma_{Q-1}$ is short, then the second sum in the definition of the coloring $c^*_X$
may change the value.
Since the distance from $x_Q$ to the final endpoint of the $\beta_i$ included in $Q$ is bounded by  $8\cdot f(r)<99\cdot f(r)$, the second sum may change the value by at most one.
\end{pf}

The following Lemma claims that the $c^*_X$-coloring of $r$-neighborhood of any piece $P$ (except for $r$-neighborhood of the set $\CC_P$) is determined by the $c^*_X$-coloring of the piece $P$.

\begin{Lem}\label{lemma color of r-close projection}
If $d(x,P)\le r$ and $\pi_P(x)\notin \CC_P$, then $c^*_X(x)=c^*_X(\pi_P(x))$.
\end{Lem}

\begin{pf}
If $\gamma$ is a geodesic from $\widehat{x}$ to $x$, then it passes through both $x_P$ and $\pi_P(x)$ (notice that $x_P\ne \pi_P(x)$).
Since the subpath $\beta$ of $\gamma$ from $\pi_P(x)$ to $x$ has length $\le r$, it cannot contain a long piece of $\gamma$, therefore
the first sum in the definition of the coloring $c^*_X$ gives the same value for $\pi_P(x)$ and for $x$.
Since the piece of $\gamma$ from $x_P$ to $\pi_P(x)$ is long, the contribution of the subpath $\beta$ to the second sum in the definition of the coloring $c^*_X$ has value $\left\lfloor \frac{|\lfloor\beta\rfloor|}{99\cdot f(r)}\right\rfloor\le\left\lfloor \frac{r}{99\cdot f(r)}\right\rfloor=0$.
\end{pf}

\begin{Lem}\label{lemma color of projected chain}
If an $r$-chain of points $\{x_i\}_{i=0}^m$ in $X$ of the same $c^*_X$-color $\ccc$
connects two points $x=x_0$ and $y=x_m$ of some piece $P$ of $X$,
then $\{\pi_P(x_i)\}_{i=0}^m$ is an $r$-chain of points of the same $c^*_X$-color $\ccc$
connecting $x$ and $y$ in $P$ provided none of the points $\pi_P(x_i)$ belongs to $\CC_P$.
\end{Lem}

\begin{pf}
By~\ref{lemma retraction to piece}, the projection $\pi_P$ does not increase distances, so the $\pi_P$-image of any $r$-chain is an $r$-chain.

We prove $c^*_X(\pi_P(x_k))=\ccc$ by induction on $k$.
Clearly, $\pi_P(x_0)=x_0$ and $c^*_X(\pi_P(x_0))=c^*_X(x_0)=\ccc$.
If $d(x_k,P)\le r$, then $c^*_X(x_k)=c^*_X(\pi_P(x_k))$ by~\ref{lemma color of r-close projection}.
If $d(x_k,P)> r$, then $\pi_P(x_k)=\pi_P(x_{k-1})$ by~\ref{lemma two points project to one}.
Thus $c^*_X(\pi_P(x_k))=c^*_X(\pi_P(x_{k-1}))=\ccc$.
\end{pf}

\begin{Lem}\label{lemma r-component in a piece}
If two points $x$ and $y$ of a piece $P$ have the same $c^*_X$-color and can be connected by an $r$-chain of that $c^*_X$-color in $X$, then $d(x,y)\le 36\cdot f(r)$.
\end{Lem}

\begin{pf}
Assume that $d(x,y) > 36\cdot f(r)$. Then at least one of the points $x$ and $y$ (say $x$)
is of distance $>18\cdot f(r)$ from the base point $x_P$ of the piece $P$.
Let $\{x_i\}_{i=0}^m$ be an $r$-chain of points in $X$ of the same $c^*_X$-color
connecting the points $x=x_0$ and $y=x_m$.
By~\ref{lemma retraction to piece}, the chain $\{\pi_P(x_i)\}_{i=0}^m$ is also an $r$-chain connecting the points $x$ and $y$ in $P$.

If none of the points $\pi_P(x_i)$ belongs to $\CC_P$,
then all points of the chain $\{\pi_P(x_i)\}_{i=0}^m$ have the same $c^*_X$-color by~\ref{lemma color of projected chain}.
So, the points $x$ and $y$ belong to one $r$-component of that color in $P$.
Therefore $d(x,y)\le 8\cdot f(r)$ by~\ref{lemma c*X almost equals cP}.

If there is a point $x_i$ with $\pi_P(x_i)\in \CC_P$, we consider the following number
$$k=\max\{l\mid \pi_P(x_i)\notin \CC_P\text{ for all } i\le l\}$$
Since $d(x_k,x_{k+1})\le r$ and $\pi_P(x_k)\ne \pi_P(x_{k+1})$, we have $d(x_k,P)\le r$ by~\ref{lemma two points project to one}.
Since $\pi_P(x_k)\notin \CC_P$,
we have $c^*_X(x_k)=c^*_X(\pi_P(x_k))$ by~\ref{lemma color of r-close projection}.
Then all points of the chain $\{\pi_P(x_i)\}_{i=0}^k$ connecting $x$ to $\pi_P(x_k)$
have the same $c^*_X$-color by~\ref{lemma color of projected chain}.
Therefore $d(x,\pi_P(x_k))\le 8\cdot f(r)$ by~\ref{lemma c*X almost equals cP}.
Thus $d(x,x_P)\le d(x,\pi_P(x_k))+d(\pi_P(x_k),\pi_P(x_{k+1}))+d(\pi_P(x_{k+1}),x_P)\le 8\cdot f(r)+r+8\cdot f(r)\le 17\cdot f(r)$
contradicting our assumption that $d(x,x_P)>18\cdot f(r)$.
\end{pf}

\begin{Lem}\label{lemma component along gamma}
Let $x$ be a point in $X$ and $\gamma$ be a geodesic from $\widehat{x}$ to $x$.
Suppose that there is an $r$-chain $\{x_i\}_{i=0}^m$ from $x=x_0$ to a point $x_m\in \gamma$
such that all the points of the chain, except possibly $x_m$, have the same $c^*_X$-color $\ccc$.
Denote by $\gamma'$ the subpath of $\gamma$ between $x$ and $x_m$.
Assume that if a long piece $\gamma_P$ of $\gamma$ intersects $\gamma'$, then either $\gamma_P\subset \gamma'$ or $\gamma_P\cap\gamma'=x_m$.
Then $d(x,x_m)\le 140\cdot f(r)$.
\end{Lem}

\begin{pf}
We consider the complement in $\gamma$ of the interiors of all long pieces of $\gamma$.
Since there are only finitely many long pieces of $\gamma$, this complement is a union of closed subpaths $\{\beta_i\}_{i=1}^k$ of $\gamma$.

{\bf Claim 1.} Any subpath $\beta'$ of any intersection $\beta_i\cap\gamma'$ of length $|\beta'|>10\cdot f(r)$ contains a point of color $\ccc$.

{\it Proof of Claim 1.}
Denote by $y$ the endpoint of $\beta'$ closest to $x_m$.
Let $\gamma_P$ be a short piece of $\gamma$ containing the point $y$.
The intersection $\gamma_P\cap\beta'$ is a subpath of $\beta'$ with endpoints $y$ and $y'$.
Since $\gamma_P$ is short, $d(y,y')\le 8\cdot f(r)$.
By~\ref{lemma chain not missing cut point}, there is an element $x_t$ of the $r$-chain
such that $d(x_t,y')\le r$ and any geodesic from $x_t$ to $x_m$ passes through $y'$.
The only reason for the point $y'$ to have $c^*_X$-color different from $c^*_X(x_t)=\ccc$ might be due to the change in the second sum of the definition of $c^*_X$-coloring as we go from $y'$ to $x_t$. Suppose this change happens.
Then the same change will happen if we go from $y'$ along $\beta'$ for the distance $d(x_t,y')\le r$
thus arriving at a point on $\beta'$ of color $\ccc$.

{\bf Claim 2.} If there is a long piece $\gamma_Q$ contained in $\gamma'$, then $|\gamma_Q|\le 17\cdot f(r)$,
it is the only long piece of $\gamma$ inside $\gamma'$,
and the point $x_m$ is $10\cdot f(r)$-close to $\gamma_Q$.

{\it Proof of Claim 2.}
Denote by $y$ the point $\pi_Q(x)$.
By~\ref{lemma chain not missing cut point}, there is an element $x_k$ of the $r$-chain such that $d(x_k,y)\le r$ and
there is an element $x_l$ of the $r$-chain such that $d(x_l,x_Q)\le r$.
By~\ref{lemma color of r-close projection}, $c^*_X(y)=c^*_X(x_k)=\ccc$.
Denote by $H$ the $r$-chain $y,x_k,x_{k+1},\dots,x_{l-1},x_l,x_Q$.
By~\ref{lemma retraction to piece}, the projection $\pi_Q(H)$ is also an $r$-chain from $y$ to $x_Q$.
Consider the first point $x_s$ of the chain $H$ such that $\pi_Q(x_s)\in \CC_Q$.
By~\ref{lemma two points project to one}, $d(x_{s-1},Q)\le r$ and by~\ref{lemma color of r-close projection}, $c^*_X(\pi_Q(x_{s-1}))=c^*_X(x_{s-1})=\ccc$.
Therefore $y,x_{k+1},\dots,x_{s-1},\pi_Q(x_{s-1})$ is an $r$-chain of color $\ccc$.
By~\ref{lemma color of projected chain}, the chain $y,\pi_Q(x_{k+1}),\dots,\pi_Q(x_{s-1}),\pi_Q(x_{s-1})$ is also an $r$-chain of color $\ccc$
(thus by~\ref{lemma c*X almost equals cP}, $d(y,\pi_Q(x_{s-1}))\le 8\cdot f(r)$).
If the set $\CCoQ$ had $c^*_X$-color $\ccc$, the chain $y,\pi_Q(x_{k+1}),\dots,\pi_Q(x_{s-1}),\pi_Q(x_{s})$ would show that the point $y$ belongs to $\CC_Q$ contradicting the definition of long piece.
So, $c^*_X(\CCoQ)\ne \ccc$.
As one goes along a geodesic from the point $\pi_Q(x_{s})$ to the point $x_s$, the $c_X^*$-color changes
from $c^*_X(\pi_Q(x_{s}))=c^*_X(\CCoQ)$ to $c^*_X(x_s)=\ccc$.
The only reason for this change is due to the change in the second sum of the definition of $c^*_X$-coloring.
Therefore $c^*_X(\CC_Q)=\ccc-1$ (mod $n+1$) and the long piece $\gamma_Q$
is preceded by some $\beta_j$ of length at least $99\cdot f(r)-d(x_s,x_Q)\ge 90\cdot f(r)$.
This would contradict Claim 1 unless $d(x_m, x_Q)<10\cdot f(r)$.
Now we estimate
$|\gamma_Q|=d(y,x_Q)\le d(y,\pi_Q(x_{s-1}))+d(\pi_Q(x_{s-1}),\pi_Q(x_{s}))+d(\pi_Q(x_{s}),x_Q)
\le 8\cdot f(r)+r+\diam(\CC_Q)\le 17\cdot f(r)$.

Now we finish the proof of lemma.
If $\gamma'$ does not contain a long piece, it is contained in some $\beta_i$ and
by Claim 1 it cannot contain subpaths of color other than $\ccc$ of length more than $10\cdot f(r)$.
By definition of the coloring $c^*_X$, the color of $\beta_i$ changes with periodicity $99\cdot f(r)$.
Thus $\gamma'$ cannot be longer than $(10+99+10)\cdot f(r)$.

If $\gamma'$ contains a long piece $\gamma_Q$, then $d(x,x_m)=d(x,\pi_Q(x))+|\gamma_Q|+d(x_Q,x_m)$.
By Claim 2, the subpath of $\gamma$ between the points $x$ and $\pi_Q(x)$ does not contain a long piece
and it changes color with periodicity $99\cdot f(r)$. By Claim 1, this subpath cannot be longer than $(10+99)\cdot f(r)$.
By Claim 2, $|\gamma_Q|\le 17\cdot f(r)$ and $d(x_Q,x_m)\le 10\cdot f(r)$.
Finally,
$$d(x,x_m)\le (10+99)\cdot f(r)+17\cdot f(r)+10\cdot f(r)< 140\cdot f(r).$$
\end{pf}

The rest of this section is devoted to the proof that $r$-components in $X$ of each $c^*_X$-color are $300\cdot f(r)$-bounded.
Consider two points $x'$ and $x''$ in $X$. Assume these point have the same $c^*_X$-color and denote it by $\ccc$.
Suppose that there exists an $r$-chain $\{x_i\}_{i=0}^l$ of color $\ccc$ from $x'=x_0$ to $x''=x_l$.
Our goal is to show that $d(x',x'')\le 300\cdot f(r)$.

If the points $x'$ and $x''$ belong to one piece, then $d(x',x'')\le 36\cdot f(r)$ by~\ref{lemma r-component in a piece}.

Suppose that the base point $\widehat{x}$ of the space $X$ and one of the points $x'$ or $x''$ (say $x''$) belong to one piece $Q$.
Denote by $z'$ the point $\pi_Q(x')$.
Consider the following number
$k=\max\{l\mid \pi_Q(x_i)=z'\text{ for all } i\le l\}$.
By~\ref{lemma two points project to one}, $d(x_k,z')\le r$.
By definition of the coloring $c^*_X$, we have
$$\ccc=c^*_X(x_k)=c_Q(\pi_Q(x_k))+\left\lfloor \frac{|\lfloor\beta\rfloor|}{99\cdot f(r)}\right\rfloor=c_Q(z')+0=c^*_X(z')$$
where $\beta$ is the part of a geodesic from $\widehat{x}$ to $x_k$ between $z'$ and $x_k$.
Therefore the chain $x'=x_0, x_1, \dots, x_k,z',x_k,x_{k+1},\dots,x_m=x''$ is an $r$-chain of color $\ccc$.
By~\ref{lemma component along gamma}, $d(x',z')\le 140\cdot f(r)$.
By~\ref{lemma r-component in a piece}, $d(z',x'')\le 36\cdot f(r)$.
So, $d(x',x'')\le d(x',z')+d(z',x'')\le 200\cdot f(r)$.

Now we assume that no two of the points $x'$, $x''$, $\widehat{x}$ belong to one piece.
Let $\gamma'$ (resp. $\gamma''$) be a geodesic from $\widehat{x}$ to $x'$ (resp. $x''$).
Let $\omega$ be a geodesic from $x'$ to $x''$.
Without loss of generality we may assume that any two of these geodesics intersect along a common subpath.
Let $x'$ and $y'$ (resp. $x''$ and $y''$) be the endpoints of the intersection $\gamma'\cap\omega$ (resp. $\gamma''\cap\omega$).
Let $\widehat{x}$ and $\widehat{y}$ be the endpoints of the intersection $\gamma'\cap\gamma''$.

Suppose $y' = y''$ and the piece of $\omega$ containing the point $y'$ is equal to $y'$.
By~\ref{lemma chain not missing cut point}, there is an element $x_k$ of the $r$-chain such that $d(x_k,y')\le r$.
Then we have $r$-chains $x'=x_0,\dots,x_k,y'$ and $y',x_k,x_{k+1},\dots,x_m=x''$.
By~\ref{lemma component along gamma}, $d(x',y')\le 140\cdot f(r)$ and $d(y',x'')\le 140\cdot f(r)$.
Thus $d(x',x'')\le 280\cdot f(r)$.

Suppose $y'\ne y''$. Then there is a simple loop formed by parts of
geodesics $\gamma'$, $\gamma''$, and $\omega$ between points $y'$, $y''$, and $\widehat{y}$.
By~\ref{def tree-graded}, this loop belongs to one piece $Q$.
We also consider here the remaining possibility that $y'=y''$ and this point
belongs to a non-trivial piece $\omega_Q$ of the geodesic $\omega$.

Denote by $z'=\pi_Q(x')$ and $z''=\pi_Q(x'')$ the endpoints of the piece $\omega_Q$.
Then $$d(x',x'') = d(x',z')+d(z',z'')+d(z'',x'')$$
By~\ref{lemma chain not missing cut point}, there is a point $x_k$ of the chain from $x'$ to $x''$ which is $r$-close to $z'$.
By~\ref{lemma component along gamma}, $d(x',z')\le 140\cdot f(r)$.
Similarly, $d(z'',x'')\le 140\cdot f(r)$.

Now we estimate $d(z',z'')$.
By~\ref{lemma chain not missing cut point}, there is a point $x_k$ of the chain such that $\pi_Q(x_k)=z'$ and $d(x_k,z')\le r$.
Similarly, there is a point $x_l$ of the chain such that $\pi_Q(x_l)=z''$ and $d(x_l,z'')\le r$.
Thus we have an $r$-chain $z',x_k, x_{k+1}, \dots, x_l,z''$. Denote this chain by $H$.
If both points $z'$ and $z''$ belong to $\CC_Q$, then $d(z',z'')\le 8\cdot f(r)$.
If the points $z'$ and $z''$ do not belong to $\CC_Q$, then $c^*_X(z')=c^*_X(x_k)=\ccc$
and $c^*_X(z'')=c^*_X(x_l)=\ccc$ by~\ref{lemma color of r-close projection};
therefore $d(z',z'')\le 36\cdot f(r)$ by~\ref{lemma r-component in a piece}.

It remains to consider the case when one of the points $z',z''$ (say $z''$) belongs to $\CC_Q$ and the other point does not belong to $\CC_Q$.
By~\ref{lemma retraction to piece}, the projection $\pi_Q(H)$ is also an $r$-chain from $z'$ to $z''$.
Consider the first point $x_s$ of the chain $H$ such that $\pi_Q(x_s)\in \CC_Q$.
By~\ref{lemma two points project to one}, $d(x_{s-1},Q)\le r$ and by~\ref{lemma color of r-close projection}, $c^*_X(\pi_Q(x_{s-1}))=c^*_X(x_{s-1})=\ccc$.
Therefore $z',x_k,x_{k+1},\dots,x_{s-1},\pi_Q(x_{s-1})$ is an $r$-chain of color $\ccc$.
By~\ref{lemma color of projected chain}, the chain $z',\pi_Q(x_k),\pi_Q(x_{k+1}),\dots,\pi_Q(x_{s-1}),\pi_Q(x_{s-1})$
is also an $r$-chain of color $\ccc$
(thus by~\ref{lemma c*X almost equals cP}, $d(z',\pi_Q(x_{s-1}))\le 8\cdot f(r)$).
Now we estimate
$d(z',z'')\le d(z',\pi_Q(x_{s-1}))+d(\pi_Q(x_{s-1}),\pi_Q(x_{s}))+d(\pi_Q(x_{s}),z'')
\le 8\cdot f(r)+r+\diam(\CC_Q)\le 17\cdot f(r)$.

Finally,
$$d(x',x'') = d(x',z')+d(z',z'')+d(z'',x'')\le 140\cdot f(r)+ 17\cdot f(r) +140\cdot f(r)<300\cdot f(r).$$

\end{document}